\documentclass[12pt]{article}
\usepackage{amsmath}
\usepackage{amsthm}

\theoremstyle{definition}
\newtheorem*{acknowledgment}{Acknowledgment}

\numberwithin{equation}{section}
\begin{document}
\def\R{{\mathbf{R}}}
\def\C{{\mathbf{C}}}
\def\N{{\mathbf{N}}}
\def\Ima{\operatorname{Im}}
\def\Rea{\operatorname{Re}}
\title{Distribution of zeros of polynomials with positive coefficients}
\author{Walter Bergweiler, Alexandre Eremenko\thanks{Supported by NSF grant
DMS-1361836.}}
\maketitle
\begin{abstract} We describe the limit zero distributions of
sequences of polynomials with positive coefficients.
\smallskip

MSC Primary: 30C15, 26C10, secondary: 31A05.

\smallskip

Keywords: polynomial, positive coefficient,
location of zeros, potential, empirical measure. 
\end{abstract}

\section{Introduction and results}
In this paper we answer the following question of Ofer Zeitouni
and Subhro Ghosh \cite{OZ}, which arises in the study of
zeros of random polynomials \cite{OZ2}.

Let $P$ be a polynomial. Consider the discrete
probability measure $\mu[P]$ in the plane which has an atom of
mass $m/\deg P$ at every zero of $P$ of multiplicity $m$.
It is called the 
``empirical measure'' in the theory
of random polynomials. 

Let $\mu_n$ be a sequence of empirical measures of some polynomials with
positive coefficients, and suppose that $\mu_n\to\mu$ weakly. 
The question is how to characterize all
possible limit measures $\mu$.
We give such a characterization in terms of logarithmic potentials.
\vspace{.1in}

\noindent
{\bf Theorem 1.} {\em
For a measure $\mu$ to be a limit of empirical measures of
polynomials with positive coefficients, it is necessary
and sufficient that the following conditions are satisfied:

$\mu$ is symmetric with respect to the complex conjugation,
$\mu(\C)\leq 1$,
and the potential
\begin{equation}\label{potential}
u(z)=\int_{|\zeta|\leq 1}\log|z-\zeta|d\mu(\zeta)+\int_{|\zeta|>1}\log
\left|1-\frac{z}{\zeta}\right|d\mu(\zeta)
\end{equation}
has the property
\begin{equation}\label{0}
u(z)\leq u(|z|).
\end{equation}}
\vspace{.1in}

The potential in Theorem 1 converges for every positive measure
with the property $\mu(\C)<\infty$ to a subharmonic
function $u\not\equiv-\infty$.
If 
$$\int_{|\zeta|>1}\log|\zeta|d\mu(\zeta)<\infty\quad\mbox{or}\quad
\int_{|\zeta|<1}\log\frac{1}{|\zeta|}d\mu(\zeta)<\infty,$$
then the definition of $u$ in Theorem 1 can be simplified to
$$\int_{\C}\log|z-\zeta|d\mu(\zeta)\quad\mbox{or}\quad
\int_{\C}\log\left|1-\frac{z}{\zeta}\right|d\mu(\zeta),$$
respectively. When these integrals exist, they differ
from the potential \eqref{potential} only by additive constants.

Obrech\-koff \cite{O} proved that empirical measures of polynomials
with non-negative coefficients satisfy
\begin{equation}\label{00}
\mu(\{ z\in\C\backslash\{0\}\colon |\arg z|\leq\alpha\})\leq\frac{2\alpha}{\pi}\mu(\C\backslash\{0\}),
\quad 0\leq\alpha\leq\pi/2.
\end{equation}
We call this the Obrech\-koff inequality. 
The limits of these measures also satisfy
(\ref{00}).

Combining our result with  Obrech\-koff's theorem we conclude that (\ref{0})
and symmetry of the measure imply (\ref{00}).
In particular we find that Obresch\-koff's inequality is satisfied not only
by polynomials with non-negative coefficients,
but more generally by polynomials
satisfying
\begin{equation}\label{ju}
|f(z)|\leq f(|z|),\quad z\in\C.
\end{equation}
The converse does not hold; that is,
the inequalities \eqref{ju} and  (\ref{0}) do not follow from
Obrech\-koff's inequality. Indeed, let 
$$P(z)=(z^2+1)^m(z^2-2z\cos\beta+1)$$
This polynomial has roots of multiplicity $m$ at $\pm i$, and simple roots at
$\exp(\pm i\beta)$. Obrech\-koff's inequality is satisfied
if $\beta\geq\pi/(2m+2)$. On the other hand, $P(1)<|P(-1)|$ for all $m$ and
$\beta\in(0,\pi/2)$.

We note that Obrech\-koff's inequality is best possible~\cite{EF}.
For other results on the roots of polynomials with positive coefficients
we refer to \cite{B}.

An important ingredient in our proof is the following
theorem of De Angelis \cite{Angelis}.
\vspace{.1in}

\noindent
{\bf Theorem A.} {\em Let 
\begin{equation}\label{poly}
f(z)=a_0+\ldots+a_dz^d,\quad a_0>0,\; a_d>0,
\end{equation}
be a real polynomial. The following conditions are equivalent:
\begin{itemize}
\item[$(i)$] There exists a positive integer $m$ such that all coefficients of
$f^m$ are strictly positive.

\item[$(ii)$] There exists a positive integer $m_0$ such that for all $m\geq m_0$,
all coefficients of $f^m$ are strictly positive.

\item[$(iii)$] The inequalities
\begin{equation}\label{1}
|f(z)|<f(|z|), \quad z\not\in [0,\infty),
\end{equation}
and
\begin{equation}\label{2}
\quad a_1>0,\quad a_{d-1}>0
\end{equation}
hold.
\end{itemize}
 }
\vspace{.1in}

\begin{acknowledgment}
We thank Ofer Zeitouni for helpful comments on this paper, as well as
John P.\ D'Angelo, David Handelman and Alan Sokal for useful discussions on Theorem~A.
\end{acknowledgment}

\section{Proof of Theorem 1}

We use some facts about subharmonic functions
and potential theory which can be found in \cite{H}.
For the reader's convenience, they are stated in the Appendix.

We recall that the Riesz measure of a subharmonic function $u$
is $(2\pi)^{-1}\Delta u$, where the Laplacian is understood
as a Schwartz distribution.
In particular the empirical measure of a polynomial $P$
of degree $d$ is
the Riesz measure of the subharmonic function $(\log|P|)/d$.
For the general properties
of convergence of subharmonic functions we refer to
\cite[Theorem 3.2.13]{H}. This result  will be used repeatedly and 
is stated for the convenience of the reader as Theorem B in the Appendix. 

The function $u$ given by (\ref{potential}) satisfies
\begin{equation}\label{Ologz}
u(z)\leq O(\log|z|),\quad z\to\infty.
\end{equation}
In turn, it is well known that every subharmonic function $u$ in the plane which
satisfies~\eqref{Ologz}
can be represented in the form (\ref{potential})
plus a constant.
We will call functions of this form simply ``potentials'';
see, for example \cite[Theorem 4.2]{Hayman2}  (case $q=0$).

\vspace{.1in}

\noindent
{\em Proof of Theorem 1.}
For a subharmonic function $u$ we put
$$B(r,u)=\max_{|z|\leq r}u(z)$$
and notice that condition (\ref{0}) can be rewritten as
\begin{equation}\label{BB}
B(r,u)=u(r),\quad r\geq 0,
\end{equation}
in view of the Maximum Principle. This implies that $u(r)$ is strictly
increasing for non-constant subharmonic functions $u$ satisfying (\ref{0}).
Moreover, the Hadamard Three Circles Theorem implies that $u(r)=B(r,v)$
is convex with respect to $\log r$, so $u(r)$ is continuous for $r>0$.

First we prove the necessity of our conditions.
Let $f_n$ be a sequence of polynomials with non-negative coefficients.
Then $u_n=\log|f_n|/\deg f_n$ are subharmonic functions whose Riesz measures
$\mu_n$ are the empirical measures of $f_n$.
As the $\mu_n$ are probability measures, every sequence contains a subsequence
for which the weak limit $\mu$ exists. This $\mu$ evidently satisfies
$\mu(\C)\leq 1$, and $\mu$ is symmetric with respect to complex conjugation.
Consider the potential $u$ defined by (\ref{potential}).
This is a subharmonic function,
$u\not\equiv-\infty$, and we have $u_n+c_n\to u$
for suitable constants $c_n$.
 
For a complete discussion of the mode of convergence here
we refer to 
the Appendix;
what we need is that $u_n(r)+c_n\to u(r)$
at every point $r>0$ and for all other points
$$\limsup_{n\to\infty}u_n(z)+c_n\leq u(|z|).$$
As the polynomials $f_n$ have non-negative coefficients, 
they satisfy (\ref{ju}),
and the $u_n$ satisfy (\ref{0}). Thus $u$
satisfies (\ref{0}). 

In the rest of this section we prove sufficiency.
We start with a measure $\mu$ such 
that the associated potential $u$ in \eqref{potential} satisfies
(\ref{0}) and 
\begin{equation}\label{symm}
u(z)=u(\overline{z}).
\end{equation}
The idea is to approximate $u$ by potentials of the form
$(\log|f_n|)/\deg f_n$,
where the $f_n$ are polynomials with real coefficients that satisfy the
assumptions of Theorem~A. Applying Theorem~A we find that
$f_n^m$ has positive coefficients for some $m$. But
$f_n^m$  has the same empirical measure as $f_n$, which is close to $\mu$.

If $u(z)=k\log|z|$, then we approximate $u$ with 
$$u_n(z)=k_n\log|z|+(1-k_n)\log|z+n|,$$
where $k_n$ is a sequence of rational numbers such that $k_n\to k,$
$0\leq k_n\leq 1$.
For the rest of the proof we assume that
$u(z)$ is not of the form $k\log|z|$.

The approximation of $u$ will be performed in several steps.
In each step we modify the function obtained on the previous step,
and starting with $u$ obtain subharmonic functions $u_1,\ldots,u_5$.
The corresponding Riesz measures will be denoted by  $\nu_1,\ldots,\nu_5$.
Each modification will preserve the asymptotic inequality~\eqref{Ologz}.

\vspace{.1in}

\noindent
1. Fix $\varepsilon>0$ and define
$$u_1(z)=\max\{ u(ze^{i\alpha})\colon |\alpha|\leq\varepsilon\}.$$
It is easy to see that $u$ is the potential of some finite measure, and that
$u_1\to u$ when $\varepsilon\to 0$. This implies that the Riesz measure of $u_1$
is close (in the weak topology) to that of $u$.

Evidently, $u_1$ satisfies \eqref{0} and \eqref{symm},
and $u_1(re^{i\theta})=u(r)$ for $|\theta|\leq\varepsilon$.
Thus $u_1(re^{i\theta})=u(r)$ 
does not
depend on $\theta$ for $|\theta|\leq\varepsilon$.
\vspace{.1in}

\noindent
2. Choose $\delta\in(0,\varepsilon)$ and consider the solution $v$
of the Dirichlet problem
in the sector
$$D=\{ z\colon |\arg z|<\delta\}$$
with boundary conditions $u_1(z)$
and satisfying $v(z)=O(\log|z|)$ as $z\to\infty$. 
To prove the existence and uniqueness of $v$, we map $D$ conformally onto
the upper half-plane,
and apply Poisson's formula to solve the Dirichlet problem.
The growth restriction near $\infty$ ensures that the solution of the
Dirichlet problem is unique.

Let $u_2$ be the result of ``sweeping out the Riesz measure'' of $u_1$ out
of the sector $D$. This means that
$$u_2(z)=\left\{\begin{array}{ll} v(z)&\; \text{for }z\in D,\\
u_1(z)&\;\mbox{otherwise}.
\end{array}\right.$$
Evidently, $u_2$ is subharmonic in the plane and satisfies \eqref{symm}.
We shall prove that $u_2$ also satisfies the strict version
of (\ref{0}), namely
\begin{equation}\label{0strict}
u_2(z)< u_2(|z|) \quad\text{for } z\notin [0,\infty).
\end{equation}
In order to do so,
we note first that $u_1$ is not harmonic in any neighborhood of
the positive ray. This follows since  $u_1(r)$ is not of the form
$u_1(r)=c\log r$ and $u_1(re^{i\theta})$ does not depend on $\theta$ 
 for $|\theta|\leq\varepsilon$. Because $u_1$ is subharmonic and 
$v$ is harmonic this implies that $v(r)>u_1(r)$ for $r>0$.
As $u_1$ satisfies~\eqref{0} we see 
that $u_2$ satisfies~\eqref{0strict} for $\delta\leq |\arg z|\leq\pi$.
In order to prove that 
$u_2$ satisfies~\eqref{0strict} also for $|\arg z|\leq\delta$,
let $G$ be the plane cut along the negative ray and define
$$\psi_\alpha(z)=z^{\alpha/\pi}\quad \text{for }z\in G,$$
with the branch of the power chosen such that 
$\psi(z)>0$ for $z>0$.
We claim that for $\alpha\in(\delta,\varepsilon)$, the function
$v_\alpha=u_2\circ\phi_\alpha$,
extended by continuity to the negative ray,
is subharmonic in the plane. Indeed, near the negative ray
this function does not depend on $\arg z$ and it is subharmonic
at all points except the negative ray, thus it is also
subharmonic in a neighborhood of the negative ray.

The limit of these subharmonic functions $v_\alpha$  as $\alpha\to\delta+0$
is the function $v_\delta$ which is thus subharmonic.
But the Riesz measure of this function $v_\delta$
is supported on the negative ray, thus
$$v_\delta(z)=\int_0^{1}\log|z+t|d\nu(t)+\int_{1+}^\infty\log\left|1+\frac{z}{t}\right|d\nu(t),$$
with some non-negative measure $\nu$. It is evident from this expression
that for every $r>0$ the function
$t\mapsto v_\delta(re^{it})$ is strictly decreasing on $[0,\pi]$.

Thus for every $r>0$, our function $t\mapsto u_2(re^{it})$
is strictly decreasing in the interval $[0,\delta]$.
This, together with the fact that $u_2$ satisfies \eqref{symm},
completes the proof that $u_2$ satisfies (\ref{0strict}).
\vspace{.1in}

\noindent
3. Now we approximate our function $u_2$ by a function $u_3$ which is
harmonic near $0$. We set
$$u_3(z)=u_2(z+\varepsilon).$$
Then $u_3$ is harmonic near the origin, and using (\ref{0strict})
and monotonicity of $u_2$ on the positive ray, we obtain
\[
u_3(z)=u_2(z+\varepsilon)<u_2(|z+\varepsilon|)\leq u_2(|z|+\varepsilon)=u_3(|z|)
\]
for $z\neq [0,\infty)$, so (\ref{0strict}) is satisfied by $u_3$.
\vspace{.1in}

\noindent
4. The subharmonic function
$u_3$ we constructed has the following properties:
\begin{itemize}
\item[a)] it satisfies (\ref{0strict}), 
\item[b)] it is harmonic near the origin,
\item[c)] it is harmonic in a neighborhood of the positive ray.
\end{itemize}
To construct a function which, in addition, is also harmonic 
near $\infty$ we consider the function
$$v(z)=u_3(1/z)+k\log|z|,$$ where $k=\nu_3(\C)$.
It is easy to see that this function
is subharmonic, if we extend it to $0$ appropriately.
Notice that $v$ satisfies (\ref{0strict}), and
it is harmonic in an angular sector containing the positive
ray (in fact in the sector $|\arg z|<\delta$). The function
$w(z)=v(z+\varepsilon)$ also satisfies (\ref{0strict}) by the same
argument that we used in Step 3 to show that $u_3$ satisfies~\eqref{0strict}.
Moreover, it is harmonic near the origin and near infinity. Thus 
the function $$u_4(z)=w(1/z)+k\log|z|$$ has all properties a), b), c) and
in addition
\begin{itemize}
\item[d)] it is harmonic in a punctured neighborhood of infinity.
\end{itemize}

\noindent
5. As $u_4$ is harmonic in a neighborhood of the origin,
it has a representation
$$u_4(z)=u_4(0)+\int\log\left|1-\frac{z}{\zeta}\right|d\nu_4(\zeta).$$
As $u_4$ satisfies~\eqref{symm}, we can write
$$u_4(x+iy)=u_4(0)+cx+O(z^2),\quad z=x+iy\to 0,$$
where
$$c=
\left.\frac{d}{dx} \left(\int \log\left|
1-\frac{x}{\zeta}\right|d\nu_4(\zeta)\right)\right|_{x=0}=
-\int \frac{\Rea \zeta }{|\zeta|^2}d\nu_4(\zeta).$$

Property (\ref{0strict}) of $u_4$ implies that $c\geq 0$.
We may achieve $c>0$ by adding to $u_4$
the potential $\varepsilon\log|1+z|$.
This procedure changes $c$ to $c+\varepsilon$. This also makes positive
the linear term in the expansion  at $\infty$.
Thus we obtain a function $u_5$, close to our original potential $u$
in the weak topology, which besides~\eqref{symm} and~\eqref{0strict}
also satisfies
\begin{align}
\label{ea}
u_5(x+iy)&=\nu_5(\C)\log|z|+b/x+O(z^{-2}),\quad z\to\infty,\\
\label{eb}
u_5(x+iy)&=u_5(0)+ax+O(z^2),\quad z=x+iy\to 0,
\end{align}
with positive constants $a$ and $b$.
\vspace{.1in}

\noindent
6. In our final step we replace the Riesz measure of $u_5$ by a 
nearby discrete probability
measure with finitely many atoms, each having rational mass.

Let $\mu$ be the Riesz measure of $u_5$.
If $\mu(\C)<1$ we change $\mu$ to a probability measure by
adding an atom sufficiently far at the negative ray.
Evidently, this procedure does not destroy our conditions~\eqref{symm}
and~\eqref{0strict},  
and we also still have (\ref{ea}) and (\ref{eb}) for certain
positive constants $a$ and $b$.

By our construction, the support of
$\mu$ is disjoint from the open set
$$H=\{ z\colon |\arg z|<\delta\}\cup\{ z\colon |z|<\delta\}\cup\{ z\colon |z|>1/\delta\},$$
and replacing $\delta$ by a smaller number if necessary we may assume that this also
holds after the atom on the negative ray was added.

Let $\mu_k$ be any sequence of symmetric discrete measures
each having finitely many
atoms of rational mass, supported outside $H$, and $\mu_k\to\mu$ weakly.
Let $w_k$ be the potential of $\mu_k$. Clearly the $w_k$ satisfy~\eqref{symm}.
We show that they also satisfy~\eqref{0strict}, provided $k$ is large.

First we consider small $|z|$, noting that the $w_k$ 
are harmonic for $|z|<\delta$.
For $z=re^{i\theta}$ with $0<r<\delta$ we thus have the expansion
\begin{equation}\label{wkseries}
w_k(z)=\sum_{n=0}^\infty a_{n,k}r^n\cos n\theta.
\end{equation}
Hence
\begin{equation}\label{fi}
\frac{\partial^2}{\partial\theta^2}w_k(z)=-a_{1,k}r\cos\theta+\Phi_k(z)
\end{equation}
with
$$\Phi_k(z)=-\sum_{n=2}^\infty a_{n,k}r^nn^2\cos n\theta.$$
As the $w_k$ are harmonic for $|z|<\delta$, the convergence to $u_5$
is locally uniformly there, and $\partial^2 w_k/\partial \theta^2$ also
converges there locally uniformly to  $\partial^2 u_5/\partial \theta^2$.
For $0<\eta<b$ and large $k$ we thus have  $a_{1,k}>\eta$  by (\ref{eb}).
Moreover, for $0<r_0<\delta$ there exists $C>0$ such that
$|w_k(z)|\leq C$ for $|z|=r_0$ and all $k$.
By Cauchy's inequalities we obtain $|a_{n,k}r_0^n|\leq C_1$ and hence 
$$|\Phi_k(z)|\leq C_2r^2\quad \text{for }r\leq r_0/2.$$
This inequality, together with (\ref{fi}) shows that $w_k$ satisfies
(\ref{0strict}) for $|z|<r_1$ with some $r_1$ independent of $k$.

The case of large $|z|$ is treated similarly, using (\ref{ea}) and
the transformation
\begin{equation}\label{transf}
u(z)\mapsto \log|z|+u(1/z),
\end{equation}
as we did before.
Thus there exists $r_2>0$ such that
$w_k$ satisfies \eqref{0strict} for $|z|>r_2$.

We finally consider the case that $r_1\leq |z|\leq r_2$.
Recall that by the first statement of Lemma~1,
$\partial^2u/\partial\theta^2$ is negative on the positive ray, so we have
a positive constant $c$ such that $(\partial^2/\partial\theta^2)u(re^{i\theta})<-c$ in some
angular sector
$$S:=\{ z\colon |\arg z|<\beta,\; r_1\leq |z|\leq r_2\}.$$
We conclude that 
$$L(r):=u(r)-u(re^{i\beta})\geq c_1>0\quad\mbox{for}\  r_1\leq r\leq r_2.$$
On the interval $[r_1,r_2]$ the convergence $w_k\to u$ is uniform,
because $u$ and $w_k$ are harmonic in $S$.
On the other hand, on the compact set 
$$K:=\{ z\colon  r_1\leq |z|\leq r_2,\; |\arg z|\geq\beta\}$$
we have $w_k(z)\leq u(z)+c_1/2$ for all sufficiently large $k$.
This follows from the general convergence properties of potentials
of weakly convergent measures 
summarized in the Appendix. 
We conclude that $w_k$ satisfies \eqref{0strict} also for $r_1\leq |z|\leq r_2$,
and hence for all $z\in\C$.

Now $w_k$ is the empirical measure of some polynomial 
$$f(z)=a_0+a_1z+\ldots+a_{d-1}z^{d-1}+a_d z^d,$$ 
and \eqref{0strict} implies that $f$ satisfies \eqref{1}.
Clearly, $a_0>0$ and $a_d>0$.
Moreover, since $a_{1,k}>0$ in \eqref{wkseries}, we see that 
$a_1>0$. The analogous expansion after the transformation~\eqref{transf}
yields that $a_{d-1}>0$.
Thus the hypotheses of Theorem~A are satisfied. 
Hence $f^m$ has positive coefficients for some $m$.
As the empirical measure of $f$ and $f^m$ coincide,
we see that $u_5$ is a limit of empirical measures of polynomials
with positive coefficients.
As we may choose $u_5$ arbitrarily close to our original potential $u$ by choosing
$\varepsilon$ sufficiently small, we see that $u$ is also a limit of
empirical measures of polynomials with positive coefficients.
This completes the proof.

\section*{Appendix: Convergence of potentials}

We frequently used various convergence properties of potentials
of weakly convergent measures which we state here for the reader's convenience.
An excellent reference for all this material is~\cite{H}.

Let $\mu_n\to\mu$ be a sequence of weakly convergent positive measures.
This means that for every continuous function $\phi$ with bounded support
$$\int\phi d\mu_n\to\int\phi d\mu,\quad n\to\infty.$$
If we restrict here to $C^\infty$-functions $\phi$ with bounded support,
we obtain convergence in the space $D'$ of Schwartz distributions.
Actually, for positive measures weak convergence is equivalent to
$D'$-convergence.

Now the sequence of subharmonic functions 
$$u_n(z)=\int_{|\zeta|\leq 1}\log|z-\zeta|d\mu_n(\zeta)+
\int_{|\zeta|>1}\log\left|1-\frac{z}{\zeta}\right|d\mu_n(\zeta)$$
converges in $D'$ to the potential of the limit measure $\mu$; that is,
we have
\begin{equation}\label{12}
\int \phi(z)u_j(z)dxdy\to \int \phi(z)u(z)dxdy
\end{equation}
for every test function~$\phi$.
For the convenience of the reader we include a standard argument showing this.

First note that with
$$
K(z,\zeta)=\left\{\begin{array}{ll}\log|z-\zeta|,&|\zeta|\leq 1,\\[1mm]
\log|1-z/\zeta|,&|\zeta|>1.
\end{array}\right.
$$
and 
$$ L(\zeta)=\int_{|z|\leq R}\phi(z)K(z,\zeta)dxdy$$
we have
$$\int \phi(z)u_j(z)dxdy = \int\phi(z)\int K(z,\zeta)d\mu_j(\zeta)dxdy
= \int L(\zeta)d\mu_j(\zeta).$$
Of course, this also holds with $u_j$ and $\mu_j$ replaced by $u$ and $\mu$.
Thus~\eqref{12} is equivalent to
\begin{equation}\label{13}
\int L(\zeta)d\mu_j(\zeta)\to \int L(\zeta)\, d\mu(\zeta).
\end{equation}
Since $|\log|1-w||\leq 2|w|$ for $|w|\leq 1/2$ we find that if $R>1$, then
$$
|K(z,\zeta)|\leq \frac{2R}{|\zeta|}\quad\mbox{for}\  |z|\leq R,\; |\zeta|>2R.
$$
Choosing $R$ such that the support of $\phi$ is contained in  $|z|\leq R$ we conclude that
$$L(\zeta)\leq \frac{C}{|\zeta|}\quad\text{for}\  |\zeta|>2R$$
with some constant $C$.

To show that~\eqref{13} holds
we choose $\varepsilon>0$ and fix $R_1>2R$ so
large that $C/R_1<\varepsilon/2$.
Now $L$ is continuous and we may write 
$L=L_1+L_2$ with continuous functions $L_1$ and $L_2$,
where $L_1$ has compact support and $L_2$ satisfies $L_2(\zeta)=0$ for 
$|\zeta|\leq R_1$ and $|L_2(\zeta)|\leq |L(\zeta)|\leq C/|\zeta|$ for  $|\zeta|>R_1$.
Then 
\begin{align*}
\left|\int L_2(\zeta)d\mu_j(\zeta)- \int L_2(\zeta)\, d\mu(\zeta)\right|
&\leq  \left|\int L_2(\zeta)d\mu_j(\zeta)\right|+\left| \int L_2(\zeta)\, d\mu(\zeta)\right|\\
&\leq  2\frac{C}{R_1}\leq \varepsilon
\end{align*}
since  $\mu_j(\C)\leq 1$ and  $\mu(\C)\leq 1$.
We also have $\int L_1(\zeta)d\mu_j(\zeta)\to \int L_1(\zeta)\, d\mu(\zeta)$
by the definition of weak convergence, which is equivalent to convergence in~$D'$.
We obtain \eqref{13} and hence~\eqref{12}.

We cite Theorem 3.2.13 from \cite{H} which says that this convergence
of potentials also holds in several other senses.
\vspace{.1in}

\noindent
{\bf Theorem B.} {\em Let $u_j\not\equiv-\infty$ be a sequence of 
subharmonic functions converging in $D'$ to the subharmonic function $u$.
Then the sequence is uniformly bounded from above on any compact set.
For every $z$ we have
\begin{equation}\label{h}
\limsup_{n\to\infty}u_n(z)\leq u(z).
\end{equation}
More generally, if $K$ is a compact set, and $f\in C(K)$, then
$$\limsup_{n\to\infty}\sup_K(u_n-f)\leq\sup_K(u-f).$$
If $d\sigma$ is a positive measure with compact support such that the potential
of $d\sigma$ is continuous, then there is equality in \eqref{h}
and $u(z)>-\infty$ for almost every $z$ with respect to $d\sigma$. Moreover, 
$u_jd\sigma\to u\,d\sigma$ weakly.}
\vspace{.1in}

In this paper we deal with subharmonic functions satisfying (\ref{BB}),
so $u(r)$ is increasing and convex with respect to $\log r$ on
$(0,\infty)$. Choosing the length element on $[0,R]$ as $d\sigma$
in Theorem B, we conclude that $u_n\to u$ almost everywhere on
the positive ray. For convex functions with respect to the logarithm
this is equivalent to the uniform convergence on compact subsets of
$(0,\infty)$. 
In particular,  $u_n(r)\to u(r)$ at every point $r>0$. As the $u_n$ 
satisfy~(\ref{0}), we conclude that 
$$\limsup_{n\to\infty}u_n(re^{i\theta})\leq u(r).$$
Choosing the uniform measure on the circle $|z|=r$ as $d\sigma$ in Theorem B,
we conclude that $u(re^{i\theta})\leq u(r)$ almost everywhere
with respect to $d\sigma$. As $u$ is upper semi-continuous, we conclude
that $u(re^{i\theta})\leq u(r)$. Thus (\ref{0}) is preserved in the limit.

{\em
\noindent
W. B.: Mathematisches Seminar, CAU Kiel,

Ludewig-Meyn-Str. 4, 24098 Kiel, Germany

bergweiler{@}math.uni-kiel.de
\vspace{.1in}

\noindent
A. E.: Department of Mathematics

Purdue University

West Lafayette, IN 47907 USA

eremenko{@}math.purdue.edu
\vspace{.1in}

}

\end{document}